\documentclass[12pt,twoside,a4paper]{amsart}
\usepackage{amssymb,hyperref}
\newcommand \C {\mathbb{C}}

\newcommand \Ordo {\mathcal{O}}
\newcommand \dbar {\overline{\partial}}

\newcommand \local {\mathcal{O}_n}
\newcommand \restr[2] {{#1}_{\big |_{#2}}}
\newcommand \diprod {\diamondsuit}
\renewcommand \epsilon \varepsilon

\DeclareMathOperator {\Hom}{Hom}
\DeclareMathOperator {\ann}{ann}
\DeclareMathOperator {\End}{End}

\DeclareMathOperator {\image}{Im}

\newtheorem{theorem}{Theorem}[section]
\newtheorem{lemma}[theorem]{Lemma}
\newtheorem{example}[theorem]{Example}

\newtheorem{proposition}[theorem]{Proposition}
\theoremstyle{remark}
\newtheorem{remark}[theorem]{Remark}
\numberwithin{equation}{section}

\title[]{A residue calculus approach to the uniform Artin-Rees lemma}
\begin{document}

\author{Jacob Sznajdman}
\address{Mathematical Sciences, Chalmers University of Technology and G\"oteborg University\\S-412 96 G\"OTEBORG\\SWEDEN}
\email{sznajdma@chalmers.se}

\renewcommand{\abstractname}{Abstract}
\begin{abstract}
  The uniform Artin-Rees lemma has been proved by C. Huneke using algebraic methods. We give a new proof
  for this result in the analytic setting, using residue calculus and a method involving complexes of
  Hermitian vector bundles. We also have to introduce a type of product
  of complexes of vector bundles, which may be applicable in the solution of other division problems
  with respect to product ideals.
\end{abstract}

\maketitle

\bibliographystyle{amsalpha}

\section{Introduction}\label{intro}
The Artin-Rees lemma is a famous result in commutative algebra from the 1950s which states the following:
\begin{theorem}[Artin-Rees]
  Let $A$ be a Noetherian ring and $M$ a finitely generated $A$-module. Given an ideal $I\subset A$ and a
  submodule $N\subset M$, there exists a number $\mu$ such that
  \begin{align*}
    I^{\mu+r}M \cap N = I^r(I^\mu M \cap N),
  \end{align*}
  for all integers $r \geq 0$.
\end{theorem}
This result was used to prove the exactness of the $I$-adic completion functor, see \cite{AMCD}.
For most applications, including the one mentioned, it suffices to know that the inclusion
\begin{align*}
  I^{\mu+r}M \cap N \subset I^r N
\end{align*}
holds.
In \cite{huneke:uniform}, Huneke showed in a general setting that the Artin-Rees lemma holds in a uniform sense,
meaning that the constant $\mu$ can be chosen independently of the ideal $I$. This is a much more delicate matter than
merely showing the existence of $\mu$ for each fixed $I$. A uniform Brian\c con-Skoda theorem is shown in the same paper.
Both theorems are proved using the same theoretical framework, namely tight closure theory.

The Brian\c con-Skoda theorem was proven in 1974 by $L^2$-methods in \cite{brianconskoda}, which in turn used
Skoda's division theorem, \cite{skodal2}. Later there have appeared
algebraic proofs (see e.g.\ \cite{lipmansathaye,lipmantessier}) and proofs that use residue calculus and integral division
formulas (see e.g.\ \cite{bgvy,andersson:bsexplicit,sznajdman:elementary}). A common feature of the latter three papers is
the use of a division formula by Berndtsson, \cite{bobformula}. A machinery, which is used in the present paper,
for constructing residue currents from
vector bundles and using the residues to study ideal membership has been developed in \cite{andersson:bullsci} and \cite{prescribed}.
So called weights, that were introduced in \cite{matsaintrep1}, are useful for obtaining integral formulas that
explicitly represent ideal memberships. These were used for example in \cite{andersson:bsexplicit} and \cite{sznajdman:elementary}.
In \cite{ass}, Huneke's uniform Brian\c con-Skoda theorem was reproved in the setting of an analytic
variety using the yoga of complexes of vector bundles and residue calculus.

Our main result is an analytic proof of the uniform Artin-Rees theorem, Theorem~\ref{AR}.
As on the algebraic side, the proof bears a similarity to the proof of the Brian\c con-Skoda theorem in \cite{ass}.

\begin{theorem}\label{AR}
  Assume that $X$ is a germ of an analytic variety at a point $x$,
  that $M$ is a finitely generated module over the local
  ring $\Ordo_{X,x}$, and that $N \subset M$ is a submodule.
  Then there exists a number $\mu$ such that for any ideal $I$ of $\Ordo_{X,x}$, the inclusion 
  \begin{align*}
    I^{\mu+r}M \cap N \subset I^rN
  \end{align*}
  holds for all integers $r \geq 0$.
\end{theorem} 

We will first observe that Theorem~\ref{AR} follows from the special case $X = \C^n$ and $M = \local^{m_0}$, where
$\local = \Ordo_{\C^n,0}$.
We know that $\Ordo_{X,x} = \local/I(X)$, so $M$ has a canonical structure as an $\local$-module
and moreover, $M = \local^{m_0} /M_0$ for some integer $m_0$ and some module $M_0$, due
to the finite generation of $M$. Take $\phi \in I^{\mu+r}M \cap N$.
If $\hat{I}$ is the inverse image of $I$ under the canonical map
$\local \to \Ordo_{X,x}$, then $\phi \in I^{\mu+r}M$ means that there exists $\hat{\phi} \in \hat{I}^{\mu+r} \local^{m_0}$
so that the image of $\hat{\phi}$ in $M$ is $\phi$. Let $\hat{N} \subset \local^{m_0}$ be the inverse image of $N$.
Then $\hat{\phi} \in \hat{N}$. We now apply Theorem~\ref{AR} in $\C^n$, which gives us that
$\hat{\phi} \in \hat{I}^r \hat{N}$. Taking images under the canonical map $\local^{m_0} \to M$, we get that $\phi \in I^r N$.

The proof of Theorem \ref{AR} will be carried out in three sections. Section~\ref{algebra} deals with the algebraic setup
which involves the definition of a certain type of product of complexes of vector bundles which we call
the $\diprod$-product. This definition is inspired by \cite{andersson:bsexplicit}.
Together with residue calculus,
this product can be used to obtain membership in products of ideals, or more generally, in tensor products of submodules
of free $\local$-modules.
In Section~\ref{currents}, we associate a residue current to any
$\diprod$-product of complexes. The problem of showing that a section $\phi$ belongs to
a tensor product of submodules, is reduced to showing that $\phi$ annihilates the product residue.
This method is applied in our proof of the Artin-Rees lemma in Section~\ref{thepf}.

\section{A complex related to tensor products of modules}\label{algebra}
We begin with some preliminaries that we shall need before we can define a product of
complexes of vector bundles. We shall also prepare for the construction of a residue
current for the product complex.

Let $X$ be a neighbourhood of $0\in\C^n$. In the sequel, all vector bundles and sheaves will be over $X$.

\subsection{Superstructures}\label{superstructures}
A superstructure is a decomposition of the sections of a sheaf or vector bundle
into parts of odd and even degrees. This generalizes the construction of the exterior algebra, and similarly,
the degree determines sign changes that occur, for example, when commuting two sections.

We will not actually carry out the details of many of the constructions we are going to mention,
as the needed arguments are similar to the ones used in the constructions of the exterior algebra and tensor products
of modules and algebras. For more details, see \cite{supermanifolds}.

Given a vector bundle $E$, a superstructure means simply a
decomposition $E=E^+ \oplus E^-$ of $E$ into an even and an odd part, that is, a $\mathbb{Z}_2$-grading of $E$.
For an element $e$ in $E^+$ or in $E^-$, its degree $\deg e$, is defined so that $\deg E^+ \equiv 0$ and $\deg E^- \equiv 1 \ (\text{mod} 2)$.
\begin{example}\label{exempel}
  For a complex of vector bundles
  \begin{align*}
    \dots \to E_2 \overset{e_2}{\to} E_1 \overset{e_1}{\to} E_0 \to 0,
  \end{align*}
the total bundle $E=\bigoplus E_i$ has a natural superstructure given by
$E^+ = \bigoplus E_{2k}$ and $E^- = \bigoplus E_{2k+1}$.
\end{example}

The endomorphism bundle $\End E$ inherits a superstructure,
such that even endomorphisms leave $E^+$ and $E^-$ invariant, whereas odd endomorphisms map $E^+$ to $E^-$
and vice versa. The sheaves of smooth forms, test forms and currents
with values in $E$, which we denote by
$\mathcal{E}(E), \mathcal{D}(E)$, and $\mathcal{D'}(E)$, respectively, all have canonical
superstructures induced from $E$. That is, the degree of a section $\alpha \otimes \omega$
is the sum of the degree of $\alpha$ as a form (or current) and the degree of $\omega$
as a section of $E$, modulo $2$.


There is some more notation to settle; we let $\mathcal{E}_X$ be the sheaf of forms with values in the
trivial bundle over $X$, and define $\mathcal{D}_X$ and $\mathcal{D'}_X$ analogously. Furthermore, $C^\infty(\cdot)$
denotes smooth sections of a sheaf or vector bundle.

An interesting object is
$\mathcal{E}(\End E) = \mathcal{E}_X \otimes C^\infty(\End E)$, which is a (sheaf of a) superalgebra, where
multiplication is given by
\begin{align*}
  (\omega_1 Y_1)(\omega_2 Y_2) = (-1)^{(\deg \omega_2)(\deg Y_1)} \omega_1 \wedge \omega_2 Y_1\circ Y_2,
\end{align*}
for $\omega_{\bullet} \in \mathcal{E}_X$ and $Y_{\bullet} \in C^\infty(\End E)$.
This algebra acts on $\mathcal{E}(E)$ by 
\begin{align}\label{action}
  (\omega \psi)(\eta e) :=
  (-1)^{(\deg \psi)(\deg \eta)} (\omega \wedge \eta) (\psi e),
\end{align}
for $\omega,\eta \in\mathcal{E}_X$, $\psi\in C^\infty(\End E)$ and $e\in C^\infty(E)$.
Similarly, $\mathcal{E}(\End E)$ acts on $\mathcal{D}'(E)$.
Moreover, currents with values in $\End E$ map $E$-valued test forms to sections of $E$. Signs are then taken into
account in a 'super' sense, that is, also similarly to \eqref{action}.

Now assume we have two complexes
$E_1$ and $E_2$ which both have superstructures. Then a superstructure on $E_1 \otimes E_2$ is induced,
so that the $\mathbb{Z}_2$-degree is simply the sum of the degrees of $E_1$ and $E_2$, that is,
\begin{align*}
  \left[E_1 \otimes E_2\right]^+ = \left(E_1^+ \otimes E_2^+\right) \oplus \left(E_1^- \otimes E_2^-\right) \\
  \left[E_1 \otimes E_2\right]^- = \left(E_1^- \otimes E_2^+\right) \oplus \left(E_1^+ \otimes E_2^-\right).
\end{align*}
Given endomorphisms $\psi_i$ of $E_i$, we define
\begin{align}\label{endo_ext}
  \psi_1 (e_1 \otimes e_2) &= \psi_1 (e_1) \otimes e_2, \\
  \psi_2 (e_1 \otimes e_2) &= (-1)^{(\deg \psi_2)(\deg e_1)}e_1 \otimes \psi_2 (e_2),\notag
\end{align}
where $e_i$ are sections of $E_i$. Thus, any endomorphism of $E_i$ induces an endomorphism of $E_1 \otimes E_2$.
It is not hard to see that an odd (even) element induces an odd (even) element.
If $\psi_i$ is a form or current valued endomorphism of $E_i$, it may be extended according to \eqref{endo_ext} nevertheless.

A final remark is that the grading of $\mathcal{E}(E_1 \otimes E_2)$ (or the endomorphism bundle),
is the sum of gradings on $E_1, E_2$ and $\mathcal{E}_X$, and similarly for currents.

\subsection{Exact complexes of hermitian vector bundles}\label{exact}
We give here a short introduction to some notions that are more thoroughly explained in \cite{prescribed}.

Assume that we are given a pointwise exact complex of hermitian vectorbundles
\begin{align*}
  \dots \to E_2 \overset{f_2}{\to} E_1 \overset{f_1}{\to} E_0 \to 0,
\end{align*}
and let $(E,f)$ be the total bundle.
Then $\nabla_E=f-\dbar$ is an operator that acts on forms and currents with values in $E$.
From $\nabla_E$ we get an operator $\nabla_{\End E}$ acting on $\End E$.
It is defined so that the `super' Leibniz rule
\begin{align}\label{superleibniz}
    \nabla_E(\alpha \omega)= (\nabla_{\End E} \alpha) \omega + (-1)^{\deg \alpha}\alpha \nabla_E \omega,
\end{align}
holds for $\omega\in\mathcal{E}(E)$ and $\alpha\in\mathcal{E}(\End E)$ or $\alpha\in\mathcal{D}'(\End E)$,
where the degree is defined by the superstructure on $E$. 

Since the map $f_1$ is surjective, we know that the equation $f_1 \psi = \phi$ is always solvable,
but it may still be difficult to find an explicit solution $\psi$. Moreover, this equation is closely related to
the equation $\nabla_E \Psi = \phi$. It is thus useful to have a endomorphism valued form $u$ such 
that $\nabla_{\End E}u = 1_{\End E}$. If $\phi$ is holomorphic, it then follows from \eqref{superleibniz}
that $\nabla_E (u\Phi) = \phi$, so $u$ gives us
a simple formula for the solution of the $\nabla_E$-equation, and thus for the original equation.
 
We will now recall a construction from \cite{prescribed} of such a form $u$.
The $\Hom(E_0,E)$-component of $u$ is written $u^0$.
In this paper we are only interested in $u^0$, so for convenience, we will drop the superscript and simply write $u$.
It satisfies the relation
\begin{align*}
  \restr{\left(\nabla_{\End E} u\right)}{E_0} = 1_{E_0}. 
\end{align*}
Note that although $u$ is has values in $\Hom(E_0,E)$, it may be that $\nabla_{\End} u$ has a component with values
in $\Hom(E_1,E)$.
Let $\sigma_k: E_{k-1} \to E_k$ be the mapping of minimal norm such that it is the inverse of $f_k$ on
the image of $f_k$, and zero on the orthogonal complement of the image. Now set $\sigma = \sum_{j\geq 1} \sigma_j$.
We then have that
\begin{align}\label{udef_sigma}
  u = \sum_{j=1}^{n+1}\left(\dbar \sigma \right)^{\wedge (j-1)} \wedge \sigma_1,
\end{align}
which is an odd form, and $u_j=\left(\dbar \sigma \right)^{\wedge (j-1)} \wedge \sigma_1$ is a $(0,j-1)$-form with values in $\Hom(E_0,E_j)$.

\subsection{The diamond product}\label{di-product}
Assume that we are given $r$ complexes of holomorphic vector bundles
\begin{align*}
  \dots \to E^k_2 \overset{f^k_2}{\to} E^k_1 \overset{f^k_1}{\to} E^k_0 \to 0,
\end{align*}
where $1 \leq k \leq r$. We bestow the bundles $E^k = \bigoplus_j E^k_j$
with superstructures as in Example~\ref{exempel}. The total map for $E^k$ is $f^k = \sum_{j} f^k_j$,
which is an $E^k$-valued endomorphism.

Our aim is to define a product of these complexes, $E^1 \diprod E^2 \diprod\cdots\diprod E^r$,
which is a new complex whose total direct sum is a subcomplex of $\bigotimes_k E^k$.

The purpose of this product is to solve membership problems in products of ideals, or more generally,
tensor products of submodules of free $\local$-modules.

We define $E^1 \diprod E^2 \diprod \dots \diprod E^r$ as the complex $(H,h)$ whose components are
\begin{align*}
  H_0 &= E^1_0 \otimes E^2_0 \otimes \dots \otimes E^r_0\\
  H_k &= \bigoplus_{\alpha_1 + \dots + \alpha_r = k-1} E^1_{1+\alpha_1} \otimes E^2_{1+\alpha_2}
  \otimes \dots \otimes E^r_{1 + \alpha_r},
\end{align*}
and whose maps $h_k: H_k \to H_{k-1}$ are
\begin{align*}
  h_1 &= f^r_1 f^{r-1}_1 \dots f^1_1, \\
  h_k &= \sum_{\substack{1 \leq s \leq r\\ j \geq 2}} \restr{{f^s_j}}{H_k}.
\end{align*}
Notice that the image of $h_1$ is simply $\image f_1 \otimes \dots \otimes \image f_r$.
In the case that $E^s_0$ are all of rank 1, we can identify $\image h_1$ with the
product ideal $\image f_1\cdot\ldots\cdot\image f_r$.
The total map $h = \sum_{k \geq 1} h_k$ can be written more concisely as
\begin{align}\label{convenient}
  h = &f^1 + f^2 + \dots + f^r -\\\notag- &f^1_1 - \dots - f^r_1 + f^r_1 f^{r-1}_1 \dots f^1_1.
\end{align}

It is straightforward to see that $(H,h)$ actually is a complex.

We note that the superstructure on $H$, which is the sum of the superstructures on each factor,
coincides with the natural superstructure in Example~\eqref{exempel} if and only if $r$ is odd.
Fortunately, we may (and will) assume that $r$ is odd 
by adding a trivial factor $0\to E \to E \to 0$ to the product, where
$E$ is any vector bundle.

Assume that for each complex $E^k$, we have a form $u^k$ , so that
\begin{align}\label{nabla_uk}
  \restr{\left(\nabla_{\End E^k}u^k\right)}{E^k_0} = 1_{E^k_0},
\end{align}
as in Section~\ref{exact}.

We then define
\begin{align}\label{udef}
  u^H = u^1 \otimes u^2 \otimes \dots \otimes u^r.
\end{align}
Since we have assumed that $r$ is odd, $u^H$ is an odd form with values in $\Hom(H_0,H)$.  

\begin{proposition}\label{Uprop}
  The form $u^H$ satisfies
  \begin{align}\label{nablageneric}
    \restr{(\nabla_{\End H} u^H)}{H_0} = 1_{H_0}.
  \end{align}
  \begin{proof}
    According to \eqref{superleibniz}, we have that $\nabla_{\End H}u^H = \nabla_H \circ u^H + u^H \circ \nabla_H$, so
    \begin{align}\label{Uprop3terms}
      \nabla_{\End H}u^H = h_1 u^H + u^H h_1 + \tilde{\nabla}_{\End H}u^H,
    \end{align}
    where $\tilde{\nabla}_{\End H}$ contains the remaining terms of $\nabla_{\End H}$. More precisely,
    $\tilde{\nabla}_{\End H}$ is the associated operator on $\End H$ obtained from
    \begin{align}\label{nablatilde}
      \tilde{\nabla}_H := \sum_{k=1}^r (f^k - f^k_1) - \dbar.      
    \end{align}
    The second term of \eqref{Uprop3terms} is zero when
    restricted to $H_0$. Due to \eqref{nabla_uk}, it follows that $f^k_1 u^k_1 = 1_{E^k_0}$; in fact,
    $f_1^k u^k_1$ is the $\Hom(E^k_0,E^k_0)$-component of $\restr{(\nabla_{\End E^k} u^k)}{E^k_0}$,
    and all other components are zero.
    Therefore
    \begin{align*}
      h_1 u^H & = f^r_1 \dots f^1_1 u^1 \otimes u^2 \otimes \dots \otimes u^r = \\
      &= f^r_1 \dots f^2_1 1_{E^1_0} \otimes u^2 \otimes \dots \otimes u^r = \\
      &=\dots = (1_{E^1_0}) \otimes \dots \otimes (1_{E^r_0}) = 1_{H_0}.
    \end{align*}
    This takes care of the first term of \eqref{Uprop3terms}, so we only need to show that
    $\restr{(\tilde{\nabla}_{\End H} u^H)}{H_0}=0$. 
    Using $\tilde{\nabla}_{\End H}u^H = \tilde{\nabla}_H\circ u^H + u^H \circ \tilde{\nabla}_H$ and the fact that $r$ is odd,
    one can check that
    \begin{align*}
      &\restr{(\tilde{\nabla}_{\End H} u^H)}{H_0} =\\
      =& \sum_1^r (-1)^{k-1}u^1 \otimes \dots \otimes u^{k-1} \otimes \restr{(\nabla_{\End E^k}-f_1^k) u^k}{E^k_0} \otimes\dots\otimes u^r.
    \end{align*}
    This is zero, because on $E^k_0$, both $\nabla_{\End E^k} u^k$ and $f_1^k u^k$ are equal to $1_{E^k_0}$.
  \end{proof}
\end{proposition}

\section{Currents associated to generically exact complexes}\label{currents}
We will recall from \cite{prescribed} how one can associate residue currents to
generically exact complexes of hermitian vector bundles. 
Most properties of these currents that we need, will be stated without proof.

We begin with a complex of hermitian vector bundles
\begin{align}\label{complex}
  \dots \to H_2 \overset{h_2}{\to} H_1 \overset{h_1}{\to} H_0 \to 0,
\end{align}
which is pointwise exact outside of a proper analytic subvariety $Z$ of $X$.
We let $H=\oplus H_j$ and $h=\oplus h_j$.
The corresponding complex of locally free sheaves is
\begin{align}\label{sheafcomplex}
  \dots \to \Ordo(H_2) \overset{h_2}{\to} \Ordo(H_1) \overset{h_1}{\to} \Ordo(H_0) \to 0.
\end{align}
Suppose we start with any subsheaf $\mathcal{J}$ of $\Ordo(H_0)$. If $H_0$ is a line bundle, then locally $\Ordo(H_0) = \local$
and $\mathcal{J}$ is simply an ideal. If we choose the complex \eqref{complex} so that $\mathcal{J}$ is the image
of the map $h_1: \Ordo (H_1) \to \Ordo (H_0)$, then $R$ encodes some, or possibly all,
information about $\mathcal{J}$. More precisely, $\ann R \subset \mathcal{J}$ and equality holds if and only if
\eqref{sheafcomplex} is exact.

Outside of $Z$, we have a form $u$ with values in $\Hom(H_0,H)$ that satisfies $\restr{(\nabla_{\End H} u)}{H_0} = 1_{H_0}$.
The definition of $u$ was given in \eqref{udef_sigma}.
There is a canonical extension of $u$
to a global current $U$, which we will now define.
Let $\chi$ be any smooth function on the real line which is identically $0$ on $(-\infty,1)$ and identically $1$
on $(2,\infty)$. For any $k$-tuple of functions $F$, set $|F|^2 = |F_1|^2 + \dots + |F_k|^2$.
We choose a non-zero $F$ that vanishes precisely on $Z$, and define
\begin{align}\label{Uext}
  U=\lim_{\epsilon\to 0}U_{\epsilon}:=\lim_{\epsilon\to 0}\chi(|F|^2/\epsilon^2)u.
\end{align}
It is non-trivial to see that this limit exists; one needs to use Hironaka's resolution of singularities.
One can also define $U$ by using analytic continuation, and this approach is taken for example in \cite{prescribed}.
These two definitions are equivalent, which can be shown by an argument similar to the proof of Lemma 2, Section 3
in \cite{bjorksamuelsson}.

We shall also use a slightly more general regularization. Instead of one tuple $F$ that vanishes on $Z$, we
take a finite number of tuples $F_1, F_2, \dots , F_r$ such that the union of their zero loci covers $Z$.
We then set
\begin{align*}
  u_{\epsilon} = \chi(|F_1|^2/\epsilon_1)\chi(|F_2|^2/\epsilon_2) \cdot\dots\cdot \chi(|F_k|^2/\epsilon_r)u,
\end{align*}
where $\epsilon$ is a tuple of positive numbers, and define
\begin{align}\label{multiU}
  U=\lim_{\epsilon_1\to 0}\lim_{\epsilon_2\to 0} \dots \lim_{\epsilon_r\to 0}u_{\epsilon}.
\end{align}

\begin{remark}\label{SEP_remark}
The two definitions of $U$ given in \eqref{Uext} and \eqref{multiU} are actually equivalent.
To show it, one can make a rather standard Hironaka argument; using a suitable desingularization
one can assume that $U$ is a meromorphic form, such that the denominator is a monomial. One
can then see that both regularizations give the same result.
As a consequence, $R$ does not depend on the regularization either.
\end{remark}

Since $\restr{(\nabla_{\End H} u)}{H_0} = 1_{H_0}$ holds outside of $Z$, we have that
\footnote{
  This notation differs from \cite{prescribed}, where $R$ has additional terms. We have only kept the term that
  acts on $H_0$.
}
\begin{align}\label{URrelation}
  \restr{(\nabla_{\End H} U)}{H_0} = 1_{H_0} - R,
\end{align}
for some current $R$ with support on $Z$. The current $R$ that appears this way is called a residue current.
The annihilator of $R$, $\ann R$, is the set of all $H_0$-valued holomorphic sections $\phi$ such that $R\phi=0$.

The following proposition gives a connection between residue calculus and module membership problems. It is proved in \cite{andersson:bullsci},
but we include it here for the reader's convenience.
\begin{proposition}\label{JcontainsAnnR}
  Let $\mathcal{J}$ be the image of the map
  \begin{align*}
    h_1: \Ordo (H_1) \to \Ordo (H_0).
  \end{align*}
  Then $\ann R \subset \mathcal{J}$.
  \begin{proof}
    Assume that $\phi \in \ann R$. By $\eqref{superleibniz}$, we get that $\nabla_{H}(U\phi)= \phi - R\phi = \phi$.
    Thus $\phi$ is $\nabla_H$-exact, since $\nabla_H \psi = \phi$, where $\psi = U\phi$. Note that $h_1(\psi) = \phi$,
    but $\psi$ is current valued, and is in general not a holomorphic solution.

    There is a decomposition $\psi=\psi_1+\psi_2+\dots$ such that $\psi_j$ is a $H_j$-valued
    $(0,j-1)$-form. We note that $\nabla_H \psi=\phi$ means that
    \begin{align*}
      h_1 \psi_1 = \phi,\quad \dbar \psi_1 = h_2 \psi_2,\quad \dots , \dbar \psi_k = 0,
    \end{align*}
    for some integer $k$. From $\dbar \psi_k = 0$, we get that locally there exists a form $\eta^k$
    such that $\psi_k = \dbar \eta^k$. This, together with $\dbar \psi_{k-1} = h_k \psi_k$
    gives that $\dbar (\psi_{k-1} - h_k \eta^k) = 0$. Solving the latter $\dbar$-equation locally, we get
    a form $\eta^{k-1}$ such that $\psi_{k-1} - h_k \eta^k = \dbar \eta^{k-1}$. Substituting this into the next equation,
    $\dbar \psi_{k-2} = h_{k-1} \psi_{k-1}$, and keeping in mind that $h_{k-1} \circ h_k = 0$, we get that
    $\dbar (\psi_{k-2} - h_{k-1} \eta^{k-1}) = 0$. Thus we get $\eta^{k-2}$ satisfying
    $\psi_{k-2} - h_{k-1} \eta^{k-1} = \dbar \eta^{k-2}$. By induction, we get that
    $\dbar ( \psi_1 - h_2 \eta^2) = 0$. The section $\hat{\psi}:= \psi_1-h_2 \eta^2$ is thus
    holomorphic and $\phi = h_1 \psi_1 = h_1 \hat{\psi}$.
  \end{proof}
\end{proposition}

Since
\begin{align}\label{Rdefprep}
  \nabla_{\End H} (\chi(|F|^2/\epsilon^2)u) = \chi(|F|^2/\epsilon^2)1_{H_0} - \dbar \chi(|F|^2/\epsilon^2)\wedge u,
\end{align}
and
\begin{align*}
  \lim_{\epsilon\to 0}\left(1 - \chi(|F|^2/\epsilon^2)\right) = 0,
\end{align*}
it follows from \eqref{URrelation} that
\begin{align}\label{Rdef}
  R=\lim_{\epsilon\to 0}\dbar \chi(|F|^2/\epsilon^2)\wedge u.
\end{align}
Note that, since the limit
in \eqref{Uext} exists, \eqref{Rdefprep} gives that the limit in \eqref{Rdef} exists too.

We will now define products of currents like $U$ and $R$. 
Assume that $U^1$ comes from a complex $E^1$ that is exact outside of an analytic set $Z^1$, and that $R^2$
comes from another complex $E^2$ that is exact outside of $Z^2$. For $i=1,2$, let $F^i$ be an analytic tuple
that vanishes precisely on $Z^i$.
We~then~define
\begin{align}\label{proddef}
  U^1 \wedge R^2 &= \lim_{\epsilon_1 \to 0} \chi(|F^1|^2/\epsilon_1^2)u^1 \wedge R^2 =\\\notag
  &=-\lim_{\epsilon_1 \to 0}\lim_{\epsilon_2 \to 0} \chi(|F^1|^2/\epsilon_1^2) \dbar\chi(|F^2|^2/\epsilon_2^2)\wedge u^1 \wedge u^2\\\label{proddef2}
  R^2 \wedge U^1 &= \lim_{\epsilon_2 \to 0} \dbar \chi(|F^2|^2/\epsilon_2^2)u^2 \wedge U^1 =\\\notag
  &=\lim_{\epsilon_2 \to 0}\lim_{\epsilon_1 \to 0} \chi(|F^1|^2/\epsilon_1^2) \dbar\chi(|F^2|^2/\epsilon_2^2)\wedge u^2 \wedge u^1.
\end{align}
The fact that these limits exist is non-trivial but can be shown using a Hironaka argument, see Proposition~4 and the following remarks in
\cite{bjorksamuelsson} or Definition~7 in \cite{samuelssonlarkang}. 
We note that these products take values in $\End (E^1 \diprod E^2)$.
Furthermore, products of more than two factors are defined analogously, that is,
each factor is regularized and limits are taken from the right to the left.


The product above depends in general on the order in which we take limits, as the following example shows,
so the product is not commutative.
\begin{example}\label{noncommexample}
  Let $z$ be a coordinate for $\C$. Then
  \begin{align*}
    \frac 1z \wedge \dbar \frac 1z = 0,
  \end{align*}
  but
  \begin{align*}
    \dbar \frac 1z \wedge \frac 1z = \dbar \frac 1{z^2} = 2\pi i \frac{\partial \delta_0}{\partial z} d\overline{z}.
  \end{align*}
\end{example}


Let $E^k$, $1\leq k \leq r$ be complexes of hermitian vector bundles that are exact outside of some sets $Z^k$,
and let $u^k$ be the associated forms on $X\setminus Z^k$ satisfying \eqref{nabla_uk}.
If $H$ is the $\diprod$-product of the complexes $E^k$, we have a form $u^H$ with values in $\Hom(H_0,H)$ with an extension $U^H$.
Again, to show that the extension exists, one has to do a Hironaka argument, and it is similar to showing
that \eqref{proddef} is well defined.
We have that $u^H$ is defined outside of the union of the sets $Z^k$, and due to Proposition~\ref{Uprop},
it satisfies \eqref{nablageneric}.
Thus, we may define the residue of $U^H$ in the same way as before, that is, as in \eqref{URrelation}.
We call this residue the product residue.
We end this section with a proposition that expresses the residue of $U^H$ in terms of the currents $U^k$ and $R^k$.

\begin{proposition}\label{prodres}
  The residue $R^H$ of $U^H$ satisfies the identity
  \begin{align}\label{Rform}
    R^H=\bigoplus_{k=1}^r (-1)^{k-1} U^1 \wedge \dots \wedge U^{k-1} \wedge R^k \wedge U^{k+1} \wedge\dots\wedge U^r.
  \end{align}
  \begin{proof}
    We will now use the multi-parameter regularization \eqref{multiU} to obtain $U^H$.
    When $\epsilon_i > 0$ for $1\leq i \leq r$, we have 
    \begin{align*}
      \restr{(\nabla_{\End H} u^H_{\epsilon})}{H_0} &= \left(\nabla_{\End H} (\prod_k \chi(|F_k|^2/\epsilon_k^2)u^H)\right)_{\bigg |_{H_0}}\\&=
      \prod_k \chi(|F_k|^2/\epsilon_k^2) 1_{H_0}-\dbar \left(\prod_k \chi(|F_k|^2/\epsilon_k^2)\right) u^H.
    \end{align*}
    As in the argument leading to \eqref{Rdef}, we thus get that
    \begin{align}\label{Rhitermeravu}
      R^H = \lim_{\epsilon_1\to 0}\lim_{\epsilon_2\to 0}\dots\lim_{\epsilon_r\to 0}\left[\dbar \left(\prod \chi(|F_k|^2/\epsilon_k^2)\right) u^H\right].
    \end{align}
    Expanding \eqref{Rhitermeravu}, we get precisely \eqref{Rform}.
  \end{proof}
\end{proposition}

\begin{remark}\label{prodresremark}
  Consider the simple case where $r=2$ and the currents $U^1$ and $U^2$ are associated to principal ideals generated by functions $f$ and $g$,
  respectively.
  By taking the complexes $E^k$ as Koszul complexes, one can check that (modulo local frames),
  $U^1 = 1/f$ and $R^1 = \dbar (1/f)$, and similarly for $E^2$.
  The proposition then just says that
  \begin{align*}
    \dbar \left(\frac 1{f} \wedge \frac 1g\right) = \dbar \frac 1{f} \wedge \frac 1{g} + \frac 1{f} \wedge \dbar \frac 1{g}.
  \end{align*}
  The reason that the sign seems to be wrong is that we have removed the frames. This makes the degree of $U^i$ even, $i=1,2$.
  By Example \ref{noncommexample}, we know that the product of currents is not commutative.
  So if we interchange the positions of $f$ and $g$, we may get a completely different decomposition of the residue.
  It is therefore important to take $\diprod$-products of complexes in an order that is suitable for the problem at hand.
\end{remark}

\section{The Artin-Rees lemma}\label{thepf}
Recall that to prove Theorem \ref{AR}, we can assume without loss of generality, that $X=\C^n$ and $M=\local^{m_0}$.
Thus $I^r M$ just consists of those tuples, all of whose entries are in $I^r$.
We fix a set of generators for $I$, say $I=(a_1,a_2,\dots,a_m)$.

We will take two complexes $E^p$ and $E^N$ so that
\begin{align}\label{two_complexes}
  \image (\Ordo(E_1^p) \to \Ordo(E_0^p)) = I^r\\\notag
  \image (\Ordo(E_1^N) \to \Ordo(E_0^N)) = N.
\end{align}
Then setting $E^{tot}=E^p \diprod E^N$, we get
\begin{align}\label{tot_complex}
  \image (\Ordo(E_1^{tot}) \to \Ordo(E_0^{tot})) = I^r N,  
\end{align}
since $I^r N$ is canonically isomorphic to $I^r \otimes N$.

Let $E^N$ be a complex that corresponds to a free resolution of the $\local$-module $M/N$,
so that $\ann R^N = N$, cf.\ the comments following \eqref{sheafcomplex}. Then $E^N$ is generically exact;
it is exact on the set where $\Ordo(E^N_0)/N = M/N$ is locally free.
For $1\leq k \leq r$ we define $E^k$ to be the Koszul complex with respect to $I$,
i.e., $E^k_1$ is a trivial vector bundle of rank $m$, and $E^k_j = \bigwedge^j E^k_1$.
We now let $E^p = E^1 \diprod E^2 \diprod \ldots \diprod E^r$.
It is straight-forward to check that \eqref{two_complexes} is satisfied under these choices.

\begin{remark}\label{exact_remark}
Using only linear algebra, one can show that $E^{tot}$ has to be exact whereever its $\diprod$-factors
$E^p$ and $E^N$ are exact. We will not prove this fact, because we did not need to use the exactness of the product
when we constructed $U^{tot}$ and $R^{tot}$ in the previous section.
\end{remark}

Subsequently, we will show that any element of $I^{\mu+r}M \cap N$ annihilates the residue current $R^{tot}$ of $E^{tot}$.
By Proposition \ref{JcontainsAnnR} and \eqref{tot_complex}, the theorem then follows.

Proposition \ref{prodres} gives that
\begin{align}\label{totres}
  R^{tot}=R^p \wedge U^N - U^p \wedge R^N,
\end{align}
where $R^p = R^{E^p}, U^p=U^{E^p}, R^N=R^{E^p}$ and $U^N=U^{E^N}$.
Assume that $\phi$ is a tuple of holomorphic functions in $I^{\mu+r}M \cap N$.
Since $\phi\in N$, we have that $R^N\phi=0$, which in turn gives that $U^p_\epsilon \wedge R^N \phi =0$, since $U^p_\epsilon$ is smooth.
Hence by \eqref{proddef},
\begin{align*}
  U^p \wedge R^N\phi = \lim_{\epsilon \to 0}U^p_\epsilon \wedge R^N \phi = 0. 
\end{align*}

A lot more work will be needed to see that also the first term of \eqref{totres} is
annihilated by $\phi$.
Let $e^k_j, 1\leq j\leq m$, be an orthonormal frame for $E^k_1$, and let $e^{k*}_j$ be the dual frame.
The maps of the Koszul complexes $\delta^k: E^k_{j+1}\to E^k_j$ are
\begin{align}\label{deltadef}
  \delta^k:= \sum_{j=1}^m a_j e^{k*}_j.
\end{align}
Outside of $Z:=Z(I)$, let $\sigma^k$ be the mapping of minimal norm such that it is the inverse of $\delta^k$ on
the image of $\delta^k$, and zero on the orthogonal complement of the image. We then have
\begin{align}\label{sigmadef}
  \sigma^k = \sum_{j=1}^m \overline{a_j} e^k_j/|a|^2,
\end{align}
where $|a|^2=\sum_{j=1}^m|a_j|^2$.
Then
\begin{align}\label{ukdef}
  u^k = \sum_{j=0}^{\min(m-1,n)} \sigma^k \wedge \left(\dbar \sigma^k \right)^{\wedge (j)}.
\end{align}
The second factor $\left(\dbar \sigma^k \right)^{\wedge (j)}$ is zero if $j>n$ or if $j>m-1$. The latter
statement follows since $\delta^k \sigma^k = 1_E^k$, so taking $\dbar$ of both sides gives that the $m$ components
of $\dbar\sigma^k$ are linearly dependent. Note that due to \eqref{sigmadef}, $u^k$ is actually explicit and
its singularity is measured precisely in terms of $|a|$.

Next, we extend the form $u^p$ (defined by \eqref{udef}) associated to $E^p$ to a global current $U^p$. Although
Proposition \ref{prodres} can be applied to obtain $R^p$, it is superfluous to regularize each form $u^j$
separately, as they are all regularized by the same function. We therefore let
\begin{align}\label{updef}
  u^p_\epsilon = \chi(|a|^2/\epsilon^2) u^p =: \chi^a_\epsilon u^p.
\end{align}

By the argument leading to \eqref{Rdef}, we get that
\begin{align}\label{Rpdef}
  R^p = \lim_{\epsilon\to 0}\dbar\chi^a_\epsilon\wedge u^p.
\end{align}

For a suitable modification $X' \overset{\pi'}{\to} X$, the current $U^N$ is a finite sum of push-forwards
of terms like $\alpha/h$, where $\alpha$ is a test form and $h$ is a monomial in some local coordinates of $X'$.
This follows from Section 2 in \cite{prescribed}.
For simplicity we will assume that there is only one such term, that is,
\begin{align*}
  U^N = \pi'_*\frac{\alpha}{h}.
\end{align*}
Then $R_\epsilon^p \wedge U^N\phi$ is the push forward of
\begin{align}\label{tokill}
  \dbar \chi^a_\epsilon \wedge u^1 \wedge u^2 \wedge \dots \wedge u^r \wedge \frac {\alpha}{h}\phi,
\end{align}
where for convenience, we have omitted to indicate any pull-backs on $\chi^a_\epsilon$, $u^j$ and $\phi$ along
the map $\pi'$. We can write \eqref{tokill} as a finite sum of terms like
\begin{align}\label{typicalterm}
  \dbar \chi^a_\epsilon \wedge \sigma^1\wedge \left(\dbar \sigma^1 \right)^{\wedge j_1}\wedge\dots\wedge
  \sigma^r\wedge \left(\dbar \sigma^r \right)^{\wedge j_r}\wedge \frac {\alpha}{h}\phi.
\end{align}
For a non-zero term, the maximal value of $\sum_{i=1}^r j_i$ is $\min(m-1,n-1)$, due to the same reason as
given for \eqref{ukdef}.
From here on we let $\alpha$ denote an arbitrary smooth form and $s_j$ be local coordinates
for $X'$ in which $h$ is a monomial.

We wish to replace \eqref{typicalterm} by similar terms where the degree of the monomial $h$ is as low as
possible. To this end we will use that
\begin{align*}
  \frac{\partial}{\partial s}\left[\frac{1}{s^k}\right] = -k \left[\frac{1}{s^{k+1}}\right], \quad k \geq 1.
\end{align*}
Furthermore, the prinicipal value current $1/h$ is a tensor product
of one-variable distributions $1/s^{k_j}_j$. Therefore \eqref{typicalterm} can be split into a sum of terms
that are derivatives of terms like
\begin{align}\label{typicalterm2}
  \partial_s^\gamma\left[ \dbar \chi^a_\epsilon \wedge \sigma^1\wedge
  \left(\dbar \sigma^1 \right)^{\wedge j_1}\wedge\dots\wedge
  \sigma^r\wedge \left(\dbar \sigma^r \right)^{\wedge j_r}\right]
  \wedge \frac {\alpha\wedge ds}{s_1 \cdot s_2 \cdot\ldots\cdot s_n}\partial_s^\beta \phi,
\end{align}
where $ds = ds_1 \wedge \ldots \wedge ds_n$, and $\beta$ and $\gamma$ are multi-indices such that $|\beta|+|\gamma|$ is at most
the order of the monomial $h$. Thus, to see that the limit of \eqref{typicalterm}, as $\epsilon$ tends to zero, is zero,
it suffices to show that the limit of \eqref{typicalterm2} is zero. 

Consider a principalization $X'' \overset{\pi}{\to} X'$ of $I$, that is, a modification such that the pull-backs $\tilde{a}_j$
of the generators $a_j$ locally
are of the form $\tilde{a}_j=a_0 a'_j$, and the tuple $a'=(a_1', \ldots , a_m')$ is non-vanishing.
The form $ds/(s_1 \cdot \ldots \cdot s_n)$ is invariant (modulo holomorphic factors) under the pull-back of $\pi$, and
\begin{align*}
  \partial_s^\beta \phi \in I^{\mu+r-|\beta|}M,
\end{align*}
so
\begin{align*}
  \pi^* \partial_s^\beta \phi \in (a_0)^{\mu+r-|\beta|}\tilde{M},
\end{align*}
where $\tilde{M}$ is isomorphic to $M=\local^{m_0}$.
Under the condition $\mu \geq \min(m,n)+|\gamma|+|\beta|$, it then follows from Lemma \ref{mainlemma} below,
that the pull-back of the form in \eqref{typicalterm2} consists of terms of the form
\begin{align}\label{endterm}
   (\chi^{\pi^*a}_\epsilon)^{(1+k)}(|a_0|/\epsilon)^{2k+2} \wedge \omega,
\end{align}
where $\omega$ is integrable and $k\geq 0$.\ The $(1+k)$ here refers to taking $1+k$ derivatives of the one-variable
function $\chi$, that is, $(\chi^{\pi^*a}_\epsilon)^{(1+k)} = \chi^{(1+k)}(|\pi^*a|^2/|\epsilon|^2)$.
Note that on the support of $(\chi^{\pi^*a}_\epsilon)^{(1+k)}$,
the quotient $|a_0|/\epsilon$ is between $1$ and $2$.
Since $(\chi^{\pi^*a}_\epsilon)^{(1+k)}$ goes to zero almost everywhere, dominated convergence shows that the limit of \eqref{endterm} is zero,
which was what we wanted to prove.

\begin{lemma}\label{mainlemma}
  For some smooth forms $\alpha_k$, one has that
  \begin{align}\label{lemmaassertion}
    \pi^* &\partial_s^\gamma\left[ \dbar \chi^a_\epsilon \wedge \sigma^1\wedge \left(\dbar \sigma^1 \right)^{\wedge j_1}\wedge\dots\wedge
      \sigma^r\wedge \left(\dbar \sigma^r \right)^{\wedge j_r}\right] = \\&
    \sum_{k=0}^{|\gamma|}\frac{(\chi^{\pi^*a}_\epsilon)^{(1+k)}(|a_0|/\epsilon)^{2k+2}\alpha_k}{\overline{a_0}a_0^{|\gamma|+\min(m,n)+r-1}}.\notag
  \end{align}
  \begin{proof}
    We will expand the left hand side by Leibniz' rule, and calculate all the terms as the $|\gamma|$ derivatives fall on various factors.
    First, there are some observations to make. The form $\dbar \sigma^j$ can be written as $\nu^j + \dbar |a|^2\wedge\alpha/|a|^4$, where
    $\alpha$ denotes an arbitrary smooth form as before, and
    \begin{align}\label{nudef}
      \nu^j = \sum_{k=1}^m \overline{\partial a_k}e^j_k/|a|^2.
    \end{align}
    Because of the factor $\dbar \chi^a_\epsilon$ in \eqref{lemmaassertion}, which is divisible by $\dbar |a|^2$,
    we can replace each occurence of $\dbar \sigma^j$ by $\nu^j$.
    We let $\overline{a}^l$ represent the product of $l$ conjugated generators $\overline{a_k}$ of $I$.
    We note that $\partial_{\tilde{s}} |a|^{-2k} = \alpha \overline{a} |a|^{-2k-2}$, for one partial derivative $\partial_{\tilde{s}}$. Thus
    \begin{align}\label{nu_derivative}
      \partial_{s}^{\gamma_1} \nu^j = \alpha \sum_{k=1}^m \overline{\partial a_k}\overline{a}^{\gamma_1}e^j_k/|a|^{2\gamma_1+2},
    \end{align}
    for an arbitrary multi-index $\gamma_1$.
    The factor $\overline{a}^{\gamma_1}$ may very well be different for each $k$, but we are only concerned with the number of conjugated
    factors of $a$. The same procedure applied to $\sigma^j$ yields
    \begin{align}\label{sigma_derivative}
      \partial_s^{\gamma_2}\sigma^j = \alpha \sum_{k=1}^m\overline{a}^{\gamma_2+1}/|a|^{2\gamma_2+2}. 
    \end{align}
    Finally, calculating $\partial_s^{\gamma_3} \dbar \chi^a_\epsilon$, we see that it is a sum of terms like
    \begin{align}\label{regform_derivative}
      \alpha \epsilon^{-2p_1-2}(\chi^a_\epsilon)^{(1+p_1)}\overline{a}^{p_1}a^{1-p_2}\overline{\partial a_k},
    \end{align}
    where $p_1+p_2 \leq |\gamma_3|$ and $p_2\leq 1$ is the number of derivatives that hit $a$. 
    The sum $p_1+p_2$ is strictly less than $|\gamma_3|$ for those terms when some derivatives fall on $\alpha$,
    and the worst case is when equality occurs.
    
    We will now expand \eqref{lemmaassertion} as promised. Let $\gamma_1, \gamma_2$ and $\gamma_3$ be the multi-index of the derivatives that hit
    ${(\nu^1)}^{j_1}\wedge\dots\wedge{(\nu^r)}^{j_r}$, $\sigma^1 \wedge \dots \wedge \sigma^r$ and $\dbar \chi^a_\epsilon$, respectively.
    Using \eqref{nu_derivative}-\eqref{regform_derivative}, we see that our typical term is
    \begin{align}\label{ourtypical}
      \pi^*\left[\alpha {(\chi^a_\epsilon)}^{(1+p_1)}\epsilon^{-2p_1-2}\overline{a}^{\left(\substack{p_1+r\\+|\gamma_1|+|\gamma_2|}\right)}a^{1-p_2}
        \overline{\partial a}^{\min(m,n)}|a|^{-2\left(\substack{\min(m,n)-1+r\\+|\gamma_1|+|\gamma_2|}\right)}\right].
    \end{align}

    The next step will be to pull everything back and keep track of the number of factors of $a,\overline{a},|a|^2$ and $\overline{\partial a}$.
    It may seem that the factors $\nu^j$ are actually worse than $\sigma^j$, but when we pull-back along the principalization,
    the factors $\overline{\partial a_k}$ will in fact be of help, and the singularity of $\nu^j$ and of $\sigma^j$ are equally severe.

    We now use that $\pi^* \overline{\partial a_k} = \overline{\partial a_0} \overline a'_k + \overline{a_0}\overline{\partial a'_k}$.
    Since any term can contain at most one factor of $\overline{\partial a_0}$, \eqref{ourtypical} can be subdivided into terms of the type
    \begin{align}\label{type1}
      \frac{{(\chi^{\pi^*a}_\epsilon)}^{(1+p_1)}(|a_0|/\epsilon)^{2p_1+2}\alpha}{a_0^{|\gamma|+\min(m,n)+r-1}},
    \end{align}
    or of the type
    \begin{align}\label{type2}
      \frac{{(\chi^{\pi^*a}_\epsilon)}^{(1+p_1)}(|a_0|/\epsilon)^{2p_1+2}\overline{\partial a_0}\alpha}{\overline{a_0}a_0^{|\gamma|+\min(m,n)+r-1}}.
    \end{align}
    Both \eqref{type1} and \eqref{type2} are of the required type, although \eqref{type1} has a slightly milder singularity, as it
    contains an extra $\overline{a_0}$ factor.
  \end{proof}
\end{lemma}



\begin{thebibliography}{Dem07}

\bibitem[A03]{matsaintrep1}
M.~Andersson, \emph{Integral representation with weights {I}}, Math. Ann. 326,
  1--18 (2003).

\bibitem[A04]{andersson:bullsci}
M.~Andersson, \emph{Residue currents and ideals of holomorphic currents}, Bull. Sci. math. 128,
 481--512 (2004).


\bibitem[A06]{andersson:bsexplicit}
\bysame \emph{{Explicit versions of the Brian\c con-Skoda theorem with
  variations}}, Michigan Math. J. {\bf 54}, no. 2, 361--373 (2006).

\bibitem[ASS08]{ass}
M.~Andersson, H.~Samuelsson, and J.~Sznajdman, \emph{{On the Brian\c con-Skoda
  theorem on a singular variety}}, Ann. Inst. Fourier. {\bf 60}, 2, 417--432 (2010). 

\bibitem[AW07]{prescribed}
M.~Andersson and E.~Wulcan, \emph{{Residue currents with prescribed annihilator ideals}}, Ann.
  Sci. \'Ecole Norm. Sup. {\bf 40}, 985--1007 (2007).

\bibitem[AW10]{anderssonwulcan:decomp}
\bysame \emph{{Decomposition of residue currents}}, J.
  reine. angew. Math. {\bf 638}, 103--118 (2010).

\bibitem[AM]{AMCD}
M.F.~Atiyah and I.G.~MacDonald, \emph{Introduction to commutative algebra},
Addison-Wesley, UK, (1999).


\bibitem[BGVY93]{bgvy}
~\\C.~Berenstein R. Gay A. Vidras~A. Yger, \emph{Residue currents and bezout
  identities}, Progress in Mathematics, 114, Birkh{\"a}user Verlag, Basel,
  1993.

\bibitem[B83]{bobformula}
B.~Berndtsson, \emph{A formula for division and interpolation}, Math. Ann. 263,
  113--160 (1983).

\bibitem[BS74]{brianconskoda}
J.~Brian{\c c}on and H. Skoda, \emph{Sur la cl\^oture int\'egrale d'un id\'eal de
  germes de fonctions holomorphes en un point de {$C^n$}}, C. R. Acad. Sci.
  Paris S\'er. A 278, 949--951 (1974).

\bibitem[BSa10]{bjorksamuelsson} \emph{Regularizations of residue currents},
J. reine angew. Math. 649, 33--54 (2010).


\bibitem[D84]{supermanifolds}
B.~DeWitt, \emph{Supermanifolds}, Cambridge University press, (1984,1992)

\bibitem[H92]{huneke:uniform}
C.~Huneke, \emph{{Uniform bounds in Noetherian rings}}, Invent. Math. {\bf
  107}, 203--223 (1992).

\bibitem[LS81]{lipmansathaye}
J.~Lipman and A.~Sathaye, \emph{{Jacobian ideals and a theorem of Brian\c
  con-Skoda}}, Michigan Math. J. {\bf 28}, no. 2, 199--222 (1981).

\bibitem[LSa10]{samuelssonlarkang}
R.~L\"ark\"ang and H.~Samuelsson \emph{Various approaches to products of residue currents},
Preprint, arXiv:1005.2056v2 (2010).

\bibitem[LT81]{lipmantessier}
J.~Lipman and B. Tessier, \emph{Pseudo-rational local rings and a theorem of
  {Brian\c con-Skoda} about integral closures of ideals}, Michigan Math. J. 28,
  97--115 (1981).

\bibitem[S72]{skodal2}
H.~Skoda, \emph{Application des techniques {$L^2$} \`a la th\'eorie des
  id\'eaux d'une alg\'ebre de fonctions holomorphes avec poids}, Ann. Sci.
  \'Ecole Norm. Sup. (4) 5, 545-579 (1972).

\bibitem[S10]{sznajdman:elementary}
J.~Sznajdman, \emph{An elementary proof of the Brian\c con-Skoda theorem},
Ann. Fac. Sci. Toulouse. Vol. 19, No. 3-4, 675--685 (2010).

\end{thebibliography}
\providecommand{\bysame}{\leavevmode\hbox to3em{\hrulefill}\thinspace}
\providecommand{\MR}{\relax\ifhmode\unskip\space\fi MR }
\providecommand{\MRhref}[2]{%
  \href{http://www.ams.org/mathscinet-getitem?mr=#1}{#2}
}
\providecommand{\href}[2]{#2}

\end{document}